\documentclass[11pt,a4paper,twoside]{article}
\usepackage{latexsym,amsmath,amsthm,amssymb,amsfonts}
\usepackage{mathrsfs,multicol}
\usepackage{amscd}
\usepackage{cite}
\usepackage[ansinew]{inputenc}

\numberwithin{equation}{section}
\newcommand{\zdhy}{\allowdisplaybreak}
\newtheorem{thm}{\indent{\bf Theorem}}[section]
\newtheorem{prop}[thm]{\indent\bf Proposition}

\newtheorem{rmk}{\indent{\bf Remark}}[section]
\newtheorem{lem}{\indent\bf Lemma}[section]

\newtheorem{coro}[thm]{\indent\bf Corollary}
\catcode`@=11

\newskip\plaincentering \plaincentering=0pt plus 1000pt minus 1000pt
\def\@plainlign{\tabskip=0pt\everycr={}}
\def\eqalignno#1{\displ@y \tabskip\plaincentering
  \halign to\displaywidth{\hfil$\@lign\displaystyle{##}$\tabskip\z@skip
    &$\@lign\displaystyle{{}##}$\hfil\tabskip\plaincentering
    &\llap{$\@lign##$}\tabskip\z@skip\crcr
    #1\crcr}}
\def\leqalignno#1{\displ@y \tabskip\plaincentering
  \halign to\displaywidth{\hfil$\@lign\displaystyle{##}$\tabskip\z@skip
    &$\@lign\displaystyle{{}##}$\hfil\tabskip\plaincentering
    &\kern-\displaywidth\rlap{$\@lign##$}\tabskip\displaywidth\crcr
    #1\crcr}}
\def\plainLet@{\relax\iffalse{\fi\let\\=\cr\iffalse}\fi}
\def\plainvspace@{\def\vspace##1{\noalign{\vskip##1}}}

\def\intic@{\mathchoice{\hskip5\p@}{\hskip4\p@}{\hskip4\p@}{\hskip4\p@}}
\def\negintic@
 {\mathchoice{\hskip-5\p@}{\hskip-4\p@}{\hskip-4\p@}{\hskip-4\p@}}
\def\intkern@{\mathchoice{\!\!\!}{\!\!}{\!\!}{\!\!}}
\def\intdots@{\mathchoice{\cdots}{{\cdotp}\mkern1.5mu
    {\cdotp}\mkern1.5mu{\cdotp}}{{\cdotp}\mkern1mu{\cdotp}\mkern1mu
      {\cdotp}}{{\cdotp}\mkern1mu{\cdotp}\mkern1mu{\cdotp}}}
\newcount\intno@
\def\iint{\intno@=\tw@\futurelet\next\ints@}
\def\iiint{\intno@=\thr@@\futurelet\next\ints@}
\def\iiiint{\intno@=4 \futurelet\next\ints@}
\def\idotsint{\intno@=\z@\futurelet\next\ints@}
\def\ints@{\findlimits@\ints@@}
\newif\iflimtoken@
\newif\iflimits@
\def\findlimits@{\limtoken@false\limits@false\ifx\next\limits
 \limtoken@true\limits@true\else\ifx\next\nolimits\limtoken@true\limits@false
    \fi\fi}
\def\multintlimits@{\intop\ifnum\intno@=\z@\intdots@
  \else\intkern@\fi
    \ifnum\intno@>\tw@\intop\intkern@\fi
     \ifnum\intno@>\thr@@\intop\intkern@\fi\intop}
\def\multint@{\int\ifnum\intno@=\z@\intdots@\else\intkern@\fi
   \ifnum\intno@>\tw@\int\intkern@\fi
    \ifnum\intno@>\thr@@\int\intkern@\fi\int}
\def\ints@@{\iflimtoken@\def\ints@@@{\iflimits@
   \negintic@\mathop{\intic@\multintlimits@}\limits\else
    \multint@\nolimits\fi\eat@}\else
     \def\ints@@@{\multint@\nolimits}\fi\ints@@@}
\def\Sb{_\bgroup\vspace@
        \baselineskip=\fontdimen10 \scriptfont\tw@
        \advance\baselineskip by \fontdimen12 \scriptfont\tw@
        \lineskip=\thr@@\fontdimen8 \scriptfont\thr@@
        \lineskiplimit=\thr@@\fontdimen8 \scriptfont\thr@@
        \Let@\vbox\bgroup\halign\bgroup \hfil$\scriptstyle
            {##}$\hfil\cr}
\def\endSb{\crcr\egroup\egroup\egroup}
\def\Sp{^\bgroup\vspace@
        \baselineskip=\fontdimen10 \scriptfont\tw@
        \advance\baselineskip by \fontdimen12 \scriptfont\tw@
        \lineskip=\thr@@\fontdimen8 \scriptfont\thr@@
        \lineskiplimit=\thr@@\fontdimen8 \scriptfont\thr@@
        \Let@\vbox\bgroup\halign\bgroup \hfil$\scriptstyle
            {##}$\hfil\cr}
\def\endSp{\crcr\egroup\egroup\egroup}
\def\Let@{\relax\iffalse{\fi\let\\=\cr\iffalse}\fi}
\def\vspace@{\def\vspace##1{\noalign{\vskip##1 }}}
\def\aligned{\,\vcenter\bgroup\plainvspace@\plainLet@\openup\jot\m@th\ialign
  \bgroup \strut\hfil$\displaystyle{##}$&$\displaystyle{{}##}$\hfil\crcr}
\def\endaligned{\crcr\egroup\egroup}
\def\matrix{\,\vcenter\bgroup\plainLet@\plainvspace@
    \normalbaselines
  \m@th\ialign\bgroup\hfil$##$\hfil&&\quad\hfil$##$\hfil\crcr
    \mathstrut\crcr\noalign{\kern-\baselineskip}}
\def\endmatrix{\crcr\mathstrut\crcr\noalign{\kern-\baselineskip}\egroup
                \egroup\,}
\newtoks\hashtoks@
\hashtoks@={#}
\def\format{\crcr\egroup\iffalse{\fi\ifnum`}=0 \fi\format@}
\def\format@#1\\{\def\preamble@{#1}%
  \def\c{\hfil$\the\hashtoks@$\hfil}%
  \def\r{\hfil$\the\hashtoks@$}%
  \def\l{$\the\hashtoks@$\hfil}%
  \setbox\z@=\hbox{\xdef\Preamble@{\preamble@}}\ifnum`{=0 \fi\iffalse}\fi
   \ialign\bgroup\span\Preamble@\crcr}

\def\cases{\left\{\,\vcenter\bgroup\plainvspace@
     \normalbaselines\openup\jot\m@th
      \plainLet@\ialign\bgroup$\displaystyle{##}$\hfil&
      \quad$\displaystyle{{}##}$\hfil\crcr
      \mathstrut\crcr\noalign{\kern-\baselineskip}}
\def\endcases{\endmatrix\right.}
\newif\iftagsleft@
\tagsleft@true
\def\TagsOnRight{\global\tagsleft@false}
\def\tag#1$${\iftagsleft@\leqno\else\eqno\fi
 \hbox{\def\pagebreak{\global\postdisplaypenalty-\@M}%
 \def\nopagebreak{\global\postdisplaypenalty\@M}\rm(#1\unskip)}%
  $$\postdisplaypenalty\z@\ignorespaces}
\interdisplaylinepenalty=\@M
\def\allowdisplaybreak{\noalign{\allowbreak}}
\def\plainallowdisplaybreak@{\def\allowdisplaybreak{\noalign{\allowbreak}}}
\def\plaindisplaybreak@{\def\displaybreak{\noalign{\break}}}
\def\align#1\endalign{\def\tag{&}\plainvspace@\plainallowdisplaybreak@
\plaindisplaybreak@
  \iftagsleft@\plainlalign@#1\endalign\else
   \plainralign@#1\endalign\fi}
\def\plainralign@#1\endalign{\displ@y\plainLet@\tabskip\plaincentering
\halign to\displaywidth
     {\hfil$\displaystyle{##}$\tabskip=\z@&$\displaystyle{{}##}$\hfil
       \tabskip=\plaincentering&\llap{\hbox{\rm(##\unskip)}}\tabskip\z@\crcr
             #1\crcr}}
\def\plainlalign@
 #1\endalign{\displ@y\plainLet@\tabskip\plaincentering\halign to \displaywidth
   {\hfil$\displaystyle{##}$\tabskip=\z@&$\displaystyle{{}##}$\hfil
   \tabskip=\plaincentering&\kern-\displaywidth
        \rlap{\hbox{\rm(##\unskip)}}\tabskip=\displaywidth\crcr
               #1\crcr}}

\def\re@#1{\par\hangindent\parindent\indent\llap{#1\enspace}\ignorespaces}
\def\qfootnote#1{\edef\@sf{\spacefactor\the\spacefactor}{}#1\@sf
      \insert\footins{\let\egroup=}\footnotesize 
      \interlinepenalty100 \let\par=\endgraf
        \leftskip=0pt \rightskip=0pt
        \splittopskip=10pt plus 1pt minus 1pt \floatingpenalty=20000
   \smallskip\re@{#1}\bgroup\strut\aftergroup{\strut\egroup}\let\next}
\catcode`\@=\active
\TagsOnRight
\begin{document}
\title{\bf Estimates for the higher order buckling eigenvalues in the unit sphere
\footnote{This research is supported by NSFC of China (No.
10671181), Project of Henan Provincial department of Sciences and
Technology (No. 092300410143), and NSF of Henan Provincial Education
department (No. 2009A110010).}}
\author{Guangyue Huang,\ Xingxiao Li
\footnote{The corresponding author. Email: xxl$@$henannu.edu.cn},\
Xuerong Qi\\
{\normalsize Department of Mathematics, Henan Normal University}
\\{\normalsize Xinxiang 453007, Henan, P.R. China} }
\date{August 12, 2009}
\maketitle
\begin{quotation}
\noindent{\bf Abstract.}~We consider the higher order buckling
eigenvalues of the following Dirichlet poly-Laplacian in the unit
sphere $(-\Delta)^p u=\Lambda (-\Delta) u$ with order $p(\geq2)$. We
obtain universal bounds on the $(k+1)$th eigenvalue in terms of the
first $k$th eigenvalues independent of the domains. In particular,
for $p=2$, our result is sharp than estimates on eigenvalues of
the buckling problem obtained by Wang and Xia in \cite{wangxiacom07}.\\
{{\bf Keywords}: eigenvalue, poly-Laplacian, buckling problem, unit sphere.} \\
{{\bf Mathematics Subject Classification}: Primary 35P15, Secondary
53C20.}

\end{quotation}

\section{Introduction}

Let $\Omega$ be a connected bounded domain in an $n$-dimensional
complete Riemannian manifold $M$.

Assume that $\lambda_i$ is the $i$th eigenvalue of the Dirichlet
poly-Laplacian with order $p$:
\begin{equation}\label{Intr1}\cases
(-\Delta)^p u=\lambda u& \ \ {\rm in}\ \Omega,\\
u=\frac{\partial u}{\partial\nu}=\cdots=\frac{\partial^{p-1}
u}{\partial\nu^{p-1}}=0&\ \ {\rm on}\
\partial\Omega,
\endcases
\end{equation}
where $\Delta$ is the Laplacian in $M$ and $\nu$ denotes the outward
unit normal vector field of $\partial\Omega$. Let
$0<\lambda_1\leq\lambda_2\leq \lambda_3\leq
\cdots\rightarrow+\infty$ denote the successive eigenvalues for
\eqref{Intr1}, where each eigenvalue is repeated according to its
multiplicity. When $p=1$, it is well known that the eigenvalue
problem \eqref{Intr1} is called a fixed membrane problem and it is
called a clamped plate problem when $p=2$. For any $p$ and
$M=\mathbb{R}^n$, Cheng-Ichikawa-Mametsuka proved in
\cite{chengqingming} the following inequality of the type of Yang:
\begin{equation}\label{Intr2}\sum_{i=1}^k(\lambda_{k+1}-\lambda_i)^2\leq\frac{4p(2p+n-2)}{n^2}
\sum_{i=1}^k(\lambda_{k+1}-\lambda_i)\lambda_i.\end{equation} In
particular, when $p=1$, the inequality \eqref{Intr2} becomes the
following inequality of Yang in \cite{yang1991}:
$$\sum_{i=1}^k(\lambda_{k+1}-\lambda_i)^2\leq\frac{4}{n}
\sum_{i=1}^k(\lambda_{k+1}-\lambda_i)\lambda_i.$$

In an
excellent paper of Cheng-Ichikawa-Mametsuka \cite{chengis09}, by
introducing functions $a_i$ and $b_i$, they considered the
eigenvalue problem \eqref{Intr1} with any order $p$ and
$M=\mathbb{S}^n(1)$. They proved that
\begin{alignat}{1}\label{Intr3}
\sum_{i=1}^k(\lambda_{k+1}-\lambda_i)^2\leq&\frac{4}{n^2}
\sum_{i=1}^k(\lambda_{k+1}-\lambda_i)\left\{\left(\lambda_i^{\frac1p}
+n\right)^p-\lambda_i\right.\\
&\left.+4\left(2^p-(p+1)\right)\lambda_i^{\frac1p}\left(\lambda_i^{\frac1p}
+n\right)^{p-2}\right\}\left(\lambda_1^{\frac1p}+\frac{n^2}{4}\right).
\end{alignat}
We remark that the inequality (2.19) in \cite{cheng05} of Cheng-Yang
and inequality (4.16) in \cite{wangxia07} of Wang-Xia are included
in the inequality \eqref{Intr3}. For the related research and
important improvement in eigenvalue problem \eqref{Intr1}, we refer
to \cite{ash,ppw2,hile1980,hook90,cheng4,wang08,huangli09,
harrell,huangchen,chengim,wucao,ash2002,ash2004} and the references
therein.

Now assume that $\Lambda_i$ is the $i$th eigenvalue of the following
Dirichlet poly-Laplacian with order $p\ (\geq2)$:
\begin{equation}\label{Intr4}\cases
(-\Delta)^p u=\Lambda (-\Delta) u& \ \ {\rm in}\ \Omega,\\
u=\frac{\partial u}{\partial\nu}=\cdots=\frac{\partial^{p-1}
u}{\partial\nu^{p-1}}=0&\ \ {\rm on}\
\partial\Omega.
\endcases
\end{equation} It is well known that this problem has a discrete
spectrum $0<\Lambda_1\leq\Lambda_2\leq \Lambda_3\leq
\cdots\rightarrow+\infty$, where each eigenvalue is repeated
according to its multiplicity. When $p=2$, the eigenvalue problem
\eqref{Intr4} is called a buckling problem. By introducing a new
method to construct nice trial functions, Cheng-Yang obtained in
\cite{chengyang06} that, for $p=2$ and $M=\mathbb{R}^n$,
\begin{equation}\label{Intr5}\sum_{i=1}^k(\Lambda_{k+1}-\Lambda_i)^2
\leq\frac{4(n+2)}{
n^2}\sum_{i=1}^k(\Lambda_{k+1}-\Lambda_i)\Lambda_i.\end{equation} As
a generalization of inequality \eqref{Intr5}, Huang-Li
\cite{huang001} considered the problem \eqref{Intr4} with any order
$p$. In fact, for $M=\mathbb{R}^n$, they proved that
\begin{equation}\label{Intr6}\sum_{i=1}^k(\Lambda_{k+1}-\Lambda_i)^2
\leq\frac{4(p-1)(n+2p-2)}{n^2}\sum_{i=1}^k(\Lambda_{k+1}
-\Lambda_i)\Lambda_i^{\frac{2p-3}{p-1}
 }.
\end{equation}
In 2007, Wang and Xia \cite{wangxiacom07} considered this problem
when $p=2$ and $M=\mathbb{S}^n(1)$. They proved that, for any
$\delta>0$,
\begin{alignat}{1}\label{Intr7}
2\sum_{i=1}^k(\Lambda_{k+1}-\Lambda_i)^2\leq&
\sum_{i=1}^k(\Lambda_{k+1}-\Lambda_i)^2\left(\delta\Lambda_i
+\frac{\delta^2(\Lambda_i-(n-2))}{4(\delta\Lambda_i+n-2)}\right)\\
&+\frac{1}{\delta}\sum_{i=1}^k(\Lambda_{k+1}-\Lambda_i)\left(\Lambda_i
+\frac{(n-2)^2}{4}\right).
\end{alignat}
We remark that the right hand side of inequality \eqref{Intr7}
depends on $\delta$. In a recent paper, by introducing a new
parameter and using Cauchy inequality, Huang-Li-Cao
\cite{huangli0000} obtain the following stronger inequality than
\eqref{Intr7} which is independent of $\delta$:
\begin{alignat}{1}\label{Intr8}
&\sum_{i=1}^k(\Lambda_{k+1}-\Lambda_i)^2
\left(2+\frac{n-2}{\Lambda_i-(n-2)}\right)\\
\leq& 2\left\{\sum_{i=1}^k(\Lambda_{k+1}-\Lambda_i)^2\left(\Lambda_i
-\frac{n-2}{\Lambda_i-(n-2)}\right)\right\}^{\frac{1}{2}}\\
&\qquad \qquad\qquad
\times\left\{\sum_{i=1}^k(\Lambda_{k+1}-\Lambda_i)
\left(\Lambda_i+\frac{(n-2)^2}{4}\right)\right\}^{\frac{1}{2}}.
\end{alignat}

Motivated by the idea used in \cite{chengis09}, we
consider in this paper the eigenvalue problem \eqref{Intr4} for any integer
$p\ (\geq2)$ when $M$ is $\mathbb{S}^n(1)$. We obtain the following
results:

{\thm\label{thm} Let $\Omega$ be a connected bounded domain in an $n$-dimensional
unit sphere $\mathbb{S}^n(1)$. Assume that $\Lambda_i$ is the $i$th
eigenvalue of the eigenvalue problem \eqref{Intr4} with $p\geq2$.
Then, we have
\begin{alignat}{1}\label{Intr9}
&\sum_{i=1}^k(\Lambda_{k+1}-\Lambda_i)^2
\left(2+\frac{n-2}{\Lambda_i-(n-2)}\right)\\
\leq&2\left\{\sum_{i=1}^k(\Lambda_{k+1}-\Lambda_i)^2
\left(f(\Lambda_i,
n)-\frac{\Lambda_i}{\Lambda_i-(n-2)}\right)\right\}^{\frac{1}{2}}\\
&\hspace{2.4cm}\times\left\{\sum_{i=1}^k(\Lambda_{k+1}-\Lambda_i)
\left(\Lambda_i+\frac{(n-2)^2}{4}\right)\right\}^{\frac{1}{2}},
\end{alignat} where
$$\align
f(\Lambda_i,
n)=&{\frac1{2(n-1)}}\left(\left(\Lambda_i^{\frac1{p-1}}+n\right)^{p-1}
-\left(\Lambda_i^{\frac1{p-1}}-n+2\right)^{p-1}\right)\\
&+\frac{n}{(n-1)}\Lambda_i^{\frac1{p-1}}\left(\Lambda_i^{\frac1{p-1}}+n\right)^{p-2}
-{\frac1{n-1}}\Lambda_i^{\frac1{p-1}}\left(\Lambda_i^{\frac1{p-1}}-n+2\right)^{p-2}\\
&+2\left(2^{p-1}-p\right)\Lambda_i^{\frac1{p-1}}
\left(\Lambda_i^{\frac1{p-1}}+n\right)^{p-3}\\
&+4(2^{p-2}-(p-1))\Lambda_i^{\frac2{p-1}}\left(\Lambda_i^{\frac1{p-1}}+n\right)^{p-4}.
\endalign$$
}

{\coro\label{coro} Under the assumptions of Theorem \ref{thm}, we have
\begin{alignat}{1}\label{Intr10}
\sum_{i=1}^k(\Lambda_{k+1}-\Lambda_i)^2\leq
\sum_{i=1}^k&(\Lambda_{k+1}-\Lambda_i)\left(f(\Lambda_i,
n)-\frac{\Lambda_i}{\Lambda_i-(n-2)}\right)\\
&\hspace{4cm}\times\left(\Lambda_i+\frac{(n-2)^2}{4}\right)
\end{alignat}
and
\begin{equation}\label{Intr11}\Lambda_{k+1}\leq
S_{k+1}+\sqrt{S_{k+1}^2-T_{k+1}},
\end{equation}
\begin{equation}\label{Intr12}\Lambda_{k+1}-\Lambda_k\leq2\sqrt{S_{k+1}^2-T_{k+1}},
\end{equation}
where $$S_{k+1}=\frac1{ k}\sum_{i=1}^k\Lambda_i+\frac{1}{
2k}\sum_{i=1}^k\left(f(\Lambda_i,
n)-\frac{\Lambda_i}{\Lambda_i-(n-2)}\right)
\left(\Lambda_i+\frac{(n-2)^2}{4}\right),$$
$$T_{k+1}=\frac{1}{k}\sum_{i=1}^k\Lambda_i^2+\frac{1}{
k}\sum_{i=1}^k\Lambda_i\left(f(\Lambda_i,
n)-\frac{\Lambda_i}{\Lambda_i-(n-2)}\right)
\left(\Lambda_i+\frac{(n-2)^2}{4}\right).$$ }

{\rmk\label{rmk} When $p=2$, we have $f(\Lambda_i, n)=\Lambda_i+1$
and
$$f(\Lambda_i, n)-\frac{\Lambda_i}{\Lambda_i-(n-2)}=\Lambda_i
-\frac{n-2}{\Lambda_i-(n-2)}.$$ Hence, for $p=2$, our inequality
\eqref{Intr9} becomes the inequality \eqref{Intr8} of Huang-Li-Cao.
Moreover, the inequality \eqref{Intr9} is sharp than the inequality
\eqref{Intr7} of Wang and Xia in \cite{wangxiacom07}.
 }

\section{Proof of the main theorem}

Let $u_i$ be the $i$th orthonormal eigenfunction of the problem
\eqref{Intr4} corresponding to the eigenvalue $\Lambda_i$, that is,
$u_i$ satisfies
\begin{equation}\label{addaddadd}\cases
(-\Delta)^p u_i=\Lambda_i(-\Delta) u_i & {\rm in}\ \Omega,\\
u_i=\frac{\partial u_i}{\partial\nu}=\cdots=\frac{\partial^{p-1}
u_i}{\partial\nu^{p-1}}=0& {\rm on}\
\partial\Omega,\\
\int_{\Omega}\langle\nabla u_i,\nabla u_j\rangle=\delta_{ij}.
\endcases
\end{equation}
Let $x_1,x_2,\ldots,x_{n+1}$ be the standard Euclidean coordinate
functions of $\mathbb{R}^{n+1}$. Then the unit sphere is defined by
$$\mathbb{S}^n(1)=\left\{ (x_1,x_2,\ldots, x_{n+1})\in\mathbb{R}^{n+1}\ ;\
\sum_{\alpha=1}^{n+1}x_\alpha^2=1\right\}.$$

Then by a rather long computation and a careful analysis, we are able to derive a sequence of inequalities which can be successfully used to
prove the following key proposition of the present paper:

{\prop\label{prop}
\begin{alignat}{1}\label{Prop}
\sum_{\alpha=1}^{n+1}\int_{\Omega}&\left(\langle\nabla x_\alpha,
\nabla u_i\rangle+x_\alpha\Delta
u_i\right)(-\Delta)^{p-2}\left(\langle\nabla
x_\alpha, \nabla u_i\rangle+x_\alpha\Delta u_i\right)\\
\leq&{\frac1{2(n-1)}}\left(\left(\Lambda_i^{\frac1{p-1}}+n\right)^{p-1}
-\left(\Lambda_i^{\frac1{p-1}}-n+2\right)^{p-1}\right)\\
&+\frac{n}{(n-1)}\Lambda_i^{\frac1{p-1}}\left(\Lambda_i^{\frac1{p-1}}+n\right)^{p-2}
-{\frac1{n-1}}\Lambda_i^{\frac1{p-1}}\left(\Lambda_i^{\frac1{p-1}}-n+2\right)^{p-2}\\
&+2\left(2^{p-1}-p\right)\Lambda_i^{\frac1{p-1}}
\left(\Lambda_i^{\frac1{p-1}}+n\right)^{p-3}\\
&+4(2^{p-2}-(p-1))\Lambda_i^{\frac2{p-1}}\left(\Lambda_i^{\frac1{p-1}}+n\right)^{p-4}.
\end{alignat}
}

We should remark that the main idea in proving Proposition \ref{prop} is similar to
that in reference \cite{chengis09}. However, here in our case, it seems a little
more complicated than in the case they considered.

For functions $f$ and $g$ defined on $\overline{\Omega}$, we define
the Dirichlet inner product $(f, g)_D$ by
$$(f, g)_D=\int_\Omega \langle\nabla f,\nabla g\rangle$$
and the Dirichlet norm of $f$ by
$$\|f\|_D=\left((f, g)_D\right)^{1/2}=\left(\int_\Omega |\nabla
f|^2\right)^{\frac{1}{2}}.$$ Define $H^2_p(\Omega)$ by
$$H^2_p(\Omega)=\{f\,:\ f, |\nabla f|,\ldots,
|\nabla^p f| \,\in L^2(\Omega)\},$$ where
$$|\nabla^p f|^2=\sum_{i_1,\cdots,i_p=1}^n
|\nabla_{i_1}\nabla_{i_2}\cdots\nabla_{i_p} f|^2.$$ Then
$H^2_p(\Omega)$ is a Hilbert space with respect to the norm
$\|\cdot\|_p$:
$$\|f\|_p=\left(\int_\Omega\left(f^2+|\nabla f|^2+\cdots+|\nabla^p f|^2\right)\right)^{\frac{1}{2}}.$$
Consider the subspace $H^2_{p,D}(\Omega)$ of $H^2_p(\Omega)$ defined
by
$$H^2_{p,D}(\Omega)=\left\{f\in H^2_p(\Omega)\,:\
f=\frac{\partial f}{\partial\nu}=\cdots=\frac{\partial^{p-1}
f}{\nu^{p-1}}=0 \ \ {\rm on}\
\partial\Omega\right\}.$$
Then the operator $(-\Delta)^p$ defines a self-adjoint operator
acting on $H^2_{p,D}(\Omega)$ for the eigenvalue problem
\eqref{Intr4} and eigenfunctions $\{u_i\}_{i=1}^\infty$ defined in
\eqref{addaddadd} form a complete orthonormal basis for the Hilbert
space $H^2_{p,D}(\Omega)$. For vector-valued functions
$$F=(f_1,f_2,\ldots,f_{n+1}),\ G=(g_1,g_2,\ldots,g_{n+1})\, :\,
\Omega\rightarrow \mathbb{R}^{n+1},$$ we define the inner product
$(F,G)$ by
$$(F,G)=\int_\Omega\langle F,G\rangle=
\int_\Omega\sum_{\alpha=1}^{n+1}f_\alpha g_\alpha.$$ The norm
of $F$ is given by
$$\|F\|=(F,F)^{\frac{1}{2}}=\left(\int_\Omega\sum_{\alpha=1}^{n+1}
f_\alpha g_\alpha\right)^{\frac{1}{2}}.$$ Let
$\mathbf{H}^2_{p-1}(\Omega)$ be the Hilbert space of vector-valued
functions given by
$$\align\mathbf{H}^2_{p-1}(\Omega)=
\Big\{F=(f_1,f_2,\ldots,f_{n+1})\,:\ &f_\alpha,|\nabla f_\alpha|,\ldots,
|\nabla^{p-1} f_\alpha| \,\in L^2(\Omega),\\
&{\rm for}\ \alpha=1, \ldots, n+1 \Big\}\endalign$$ with norm
$$\|F\|_{p-1}=\left\{\|F\|^2+\int_\Omega
\left(\sum_{\alpha=1}^{n+1}|\nabla
f_\alpha|^2+\cdots+\sum_{\alpha=1}^{n+1}|\nabla^{p-1}
f_\alpha|^2\right)\right\}^{\frac{1}{2}}.$$ Observe that a vector
field on $\Omega$ can be regarded as a vector-valued function from
$\Omega$ to $\mathbb{R}^{n+1}$. Let
$\mathbf{H}^2_{p-1,D}(\Omega)\subset\mathbf{H}^2_{p-1}(\Omega)$ be a
subspace of $\mathbf{H}^2_{p-1}(\Omega)$ spanned by the
vector-valued functions  $\{\nabla u_i\}_{i=1}^\infty$ which form a
complete orthonormal basis of $\mathbf{H}^2_{p-1,D}(\Omega)$. For
any $f\in H^2_{p,D}(\Omega)$, we have $\nabla
f\in\mathbf{H}^2_{p-1,D}(\Omega)$ and for any $X\in
\mathbf{H}^2_{p-1,D}(\Omega)$, there exists a function $f\in
H^2_{p,D}(\Omega)$ such that $X=\nabla f$.

Let $x_1,x_2,\ldots,x_{n+1}$ be the standard Euclidean coordinate
functions of $\mathbb{R}^{n+1}$, and $u_i$ be the $i$-th orthonormal
eigenfunction of the problem \eqref{Intr4} corresponding to the
eigenvalue $\Lambda_i$ (see \eqref{addaddadd}). For any $\alpha=1,2,\ldots, n+1$
and each $i=1,\ldots,k$, we decompose the vector-valued functions
$x_\alpha\nabla u_i$ as
\begin{equation}\label{2sec1}x_\alpha\nabla u_i=\nabla h_{\alpha i}+W_{\alpha i},\end{equation}
where $h_{\alpha i}\in H^2_{p,D}(\Omega)$, $\nabla h_{\alpha i}$ is
the projection of $x_\alpha\nabla u_i$ in
$\mathbf{H}^2_{p-1,D}(\Omega)$, $W_{\alpha i}\perp
\mathbf{H}^2_{p-1,D}(\Omega)$. Thus we have
\begin{equation}\label{2sec2}
(W_{\alpha i},\nabla u)=\int_\Omega \langle W_{\alpha
i},\nabla u\rangle=0,\ \ {\rm for\ any} \ u\in H^2_{p,D}(\Omega).
\end{equation}
By the denseness of $H^2_{p,D}(\Omega)$ in $L^2(\Omega)$ and
$C^1(\Omega)$ is dense in $L^2(\Omega)$, we conclude that
\begin{equation}\label{2sec3}
(W_{\alpha i},\nabla h)=0,\ \ \ \forall \ h\in C^1(\Omega)\cap
L^2(\Omega),
\end{equation}
which implies from the divergence theorem that
$$
\int_\Omega h\,\, {\rm div}(W_{\alpha i})=0,
$$
where ${\rm div}(Z)$ denotes the divergence of $Z$. Consequently, we
get
\begin{equation}\label{2sec4}
{\rm div}(W_{\alpha i})=0.
\end{equation}

Define $\phi_{\alpha i}$ by
\begin{equation}\label{2sec5} \phi_{\alpha i}=h_{\alpha i}-\sum_{j=1}^kb_{\alpha ij}u_j,
\end{equation}
where $$b_{\alpha ij}=\int_\Omega x_\alpha \langle\nabla
u_i,\nabla u_j\rangle=b_{\alpha ji}.$$ Then we have $$\phi_{\alpha
i}=\frac{\partial \phi_{\alpha i}}{\partial\nu}=\cdots=\frac{\partial^{p-1}
\phi_{\alpha i}}{\partial\nu^{p-1}}=0$$ and
\begin{equation}\label{2sec6}
(\phi_{\alpha i},u_j)_D=\int_\Omega\langle\nabla\phi_{\alpha
i},\nabla u_j\rangle =0, \  \ \ {\rm  for \ any}\ j=1,\ldots,
k.\end{equation} It follows from the Rayleigh-Ritz inequality that
\begin{equation}\label{2sec7}\Lambda_{k+1}\leq
\frac{\int_\Omega \phi_{\alpha i}(-\Delta)^p\phi_{\alpha
i}}{\|\nabla \phi_{\alpha i}\|^2},\end{equation} where
$\|f\|^2=\int_\Omega |f|^2$. It is easy to see from
\eqref{2sec5} and \eqref{2sec6} that
\begin{alignat}{1}\label{2sec8}
\int_{\Omega}\phi_{\alpha i}(-\Delta)^p\phi_{\alpha i}
=&\int_{\Omega}\phi_{\alpha i}\left((-\Delta)^ph_{\alpha i}
-\sum_{j=1}^k b_{\alpha ij}\Lambda_j (-\Delta)u_j\right)\\
=&\int_{\Omega}\phi_{\alpha i}(-\Delta)^ph_{\alpha i}\\
=&\int_{\Omega}\left(h_{\alpha i}
-\sum_{j=1}^kb_{\alpha ij}u_j\right)(-\Delta)^ph_{\alpha i}\\
=&\int_{\Omega}h_{\alpha i}(-\Delta)^ph_{\alpha i}
-\sum_{j=1}^kb_{\alpha ij}\int_{\Omega}u_j(-\Delta)^ph_{\alpha i}\\
=&\int_{\Omega}h_{\alpha i}(-\Delta)^ph_{\alpha i}
-\sum_{j=1}^kb_{\alpha ij}\int_{\Omega}h_{\alpha i}(-\Delta)^pu_j\\
=&\int_{\Omega}h_{\alpha i}(-\Delta)^ph_{\alpha
i}-\sum_{j=1}^k\Lambda_jb_{\alpha ij}^2. \end{alignat} Since
\begin{equation}\label{2sec9}
\|x_\alpha\nabla u_i\|^2=\int_{\Omega}x_\alpha^2|\nabla
u_i|^2=\|\nabla h_{\alpha i}\|^2+\|W_{\alpha i}\|^2,
\end{equation}
\begin{equation}\label{2sec10}
\|\nabla h_{\alpha i}\|^2=\|\nabla \Phi_{\alpha
i}\|^2+\sum_{j=1}^kb_{\alpha ij}^2.
\end{equation}
Therefore, \eqref{2sec8} can be written as
\begin{alignat}{1}\label{2sec11}
\int_{\Omega}\phi_{\alpha i}(-\Delta)^p\phi_{\alpha i}=&
\int_{\Omega}h_{\alpha i}(-\Delta)^ph_{\alpha
i}-\Lambda_i\|x_\alpha\nabla u_i\|^2\\
&+\Lambda_i\left(\|\nabla \phi_{\alpha i}\|^2+\|W_{\alpha
i}\|^2+\sum_{j=1}^kb_{\alpha
ij}^2\right)-\sum_{j=1}^k\Lambda_jb_{\alpha ij}^2.
\end{alignat}
Inserting \eqref{2sec11} into \eqref{2sec7} yields
\begin{alignat}{1}\label{2sec12}
(\Lambda_{k+1}-\Lambda_i)\|\nabla \Phi_{\alpha i}\|^2\leq&
\int_{\Omega}h_{\alpha i}(-\Delta)^ph_{\alpha
i}-\Lambda_i\|x_\alpha\nabla u_i\|^2+\Lambda_i\|W_{\alpha
i}\|^2\\&+\sum_{j=1}^k(\Lambda_i-\Lambda_j)b_{\alpha ij}^2\\
=&p_{\alpha i}+\|\langle\nabla x_\alpha,\nabla
u_i\rangle\|^2+\Lambda_i\|W_{\alpha
i}\|^2\\&+\sum_{j=1}^k(\Lambda_i-\Lambda_j)b_{\alpha ij}^2,
\end{alignat}
where $$p_{\alpha i}=\int_{\Omega}h_{\alpha
i}(-\Delta)^ph_{\alpha i}-\Lambda_i\|x_\alpha\nabla
u_i\|^2-\|\langle\nabla x_\alpha,\nabla u_i\rangle\|^2.$$

{\lem{\em\cite{wangxia07}}
Let $$c_{\alpha ij}=\int_{\Omega}\langle Z_{\alpha
i},u_j\rangle,$$ where $Z_{\alpha i}=\nabla \langle x_\alpha,\nabla
u_i\rangle-\frac{n-2}{2}x_\alpha\nabla u_i$. Then we have $$c_{\alpha
ij}=-c_{\alpha ji}.$$ }

Note that
\begin{alignat}{1}\label{2sec13}
&-2\int_{\Omega}\langle x_\alpha\nabla u_i,Z_{\alpha
i}\rangle\\
=&-2\int_{\Omega}\langle x_\alpha\nabla u_i,\nabla
\langle x_\alpha,\nabla
u_i\rangle\rangle+(n-2)\int_{\Omega}x_\alpha^2|\nabla
u_i|^2\\
=&2\int_{\Omega}\langle x_\alpha,\nabla
u_i\rangle^2+\int_{\Omega}\langle\nabla x_\alpha^2,\nabla
u_i\rangle\Delta u_i+(n-2)\int_{\Omega}x_\alpha^2|\nabla
u_i|^2.
\end{alignat}
On the other hand, from \eqref{2sec1}, \eqref{2sec3} and
\eqref{2sec6}, we obtian
\begin{alignat}{1}\label{2sec14}
-2\int_{\Omega}\langle x_\alpha\nabla u_i,Z_{\alpha
i}\rangle
=&-2\int_{\Omega}\langle\nabla h_{\alpha
i}+W_{\alpha i},Z_{\alpha i}\rangle\\ \zdhy
=&-2\int_{\Omega}\langle\nabla h_{\alpha i},Z_{\alpha
i}\rangle+(n-2)\int_{\Omega}
\langle W_{\alpha i},x_\alpha\nabla u_i\rangle\\ \zdhy
=&-2\int_{\Omega}\langle\nabla \phi_{\alpha
i}+\sum_{j=1}^kb_{\alpha ij}\nabla u_j,Z_{\alpha
i}\rangle+(n-2)\int_{\Omega}
\langle W_{\alpha i},x_\alpha\nabla u_i\rangle\\ \zdhy
=&-2\int_{\Omega}\langle\nabla \phi_{\alpha i},Z_{\alpha
i}\rangle-2\sum_{j=1}^kb_{\alpha ij}c_{\alpha
ij}+(n-2)\| W_{\alpha i}\|^2\\ \zdhy
=&-2\int_{\Omega}\langle\nabla \phi_{\alpha i},Z_{\alpha
i}-\sum_{j=1}^kc_{\alpha ij}\nabla u_j\rangle-2\sum_{j=1}^kb_{\alpha
ij}c_{\alpha ij}\\ \zdhy
&+(n-2)\| W_{\alpha i}\|^2.
\end{alignat}
From \eqref{2sec13} and \eqref{2sec14}, we obtain
\begin{equation}\label{2sec15}
r_{\alpha i}+2\sum_{j=1}^kb_{\alpha ij}c_{\alpha
ij}=-2\int_{\Omega}\langle\nabla \phi_{\alpha i},Z_{\alpha
i}-\sum_{j=1}^kc_{\alpha ij}\nabla u_j\rangle+(n-2)\| W_{\alpha
i}\|^2,
\end{equation}
where $$r_{\alpha i}=2\int_{\Omega}\langle x_\alpha,\nabla
u_i\rangle^2+\int_{\Omega}\langle\nabla x_\alpha^2,\nabla
u_i\rangle\Delta u_i+(n-2)\int_{\Omega}x_\alpha^2|\nabla
u_i|^2.$$ Multiplying \eqref{2sec15} by
$(\Lambda_{k+1}-\Lambda_i)^2$, one obtains from the Schwarz
inequality and \eqref{2sec12} that
\begin{alignat}{1}\label{2sec16}
&(\Lambda_{k+1}-\Lambda_i)^2\left(r_{\alpha
i}+2\sum_{j=1}^kb_{\alpha ij}c_{\alpha ij}\right)\\ \zdhy
=&(\Lambda_{k+1}-\Lambda_i)^2\left(-2\int_{\Omega}\left\langle\nabla
\phi_{\alpha i},Z_{\alpha i}-\sum_{j=1}^kc_{\alpha ij}\nabla
u_j\right\rangle+(n-2)\| W_{\alpha i}\|^2\right)\\ \zdhy
\leq&\delta(\Lambda_{k+1}-\Lambda_i)^3\|\nabla \phi_{\alpha
i}\|^2+\frac{1}{\delta}(\Lambda_{k+1}-\Lambda_i)\left\|Z_{\alpha
i}-\sum_{j=1}^kc_{\alpha ij}\nabla
u_j\right\|^2\\ \zdhy
&+(n-2)(\Lambda_{k+1}-\Lambda_i)^2\| W_{\alpha i}\|^2\\ \zdhy
\leq&\delta(\Lambda_{k+1}-\Lambda_i)^2\left(p_{\alpha
i}+\|\langle\nabla x_\alpha,\nabla
u_i\rangle\|^2+\Lambda_i\|W_{\alpha
i}\|^2+\sum_{j=1}^k(\Lambda_i-\Lambda_j)b_{\alpha ij}^2\right)\\ \zdhy
&+\frac{1}{\delta}(\Lambda_{k+1}-\Lambda_i)\left(\|Z_{\alpha
i}\|^2-\sum_{j=1}^kc_{\alpha
ij}^2\right)+(n-2)(\Lambda_{k+1}-\Lambda_i)^2\| W_{\alpha i}\|^2.
\end{alignat}
Since $b_{\alpha ij}=b_{\alpha ji}$ and $c_{\alpha ij}=-c_{\alpha
ji}$, summing over $i$ from 1 to $k$ for \eqref{2sec16} yields
\begin{alignat}{1}\label{2sec17}
&\sum_{i=1}^k(\Lambda_{k+1}-\Lambda_i)^2r_{\alpha
i}\\
\leq&\sum_{i=1}^k(\Lambda_{k+1}-\Lambda_i)^2\Big(\delta p_{\alpha
i}+\delta\|\langle\nabla x_\alpha,\nabla
u_i\rangle\|^2+(\delta\Lambda_i+n-2)\|W_{\alpha
i}\|^2\Big)\\
&+\frac{1}{\delta}\sum_{i=1}^k(\Lambda_{k+1}-\Lambda_i)\|Z_{\alpha
i}\|^2.
\end{alignat}
Let $\rho$ be a positive constant. Then we have
\begin{alignat}{1}\label{2sec18}
\rho\|\langle\nabla x_\alpha,\nabla
u_i\rangle\|^2=&\rho\int_\Omega\langle\nabla x_\alpha,\nabla
u_i\rangle^2\\
=&-\rho\int_\Omega  x_\alpha {\rm div}(\langle\nabla
x_\alpha,\nabla u_i\rangle\nabla u_i)\\
=&-\rho\int_\Omega\langle x_\alpha\nabla u_i,\nabla \langle\nabla
x_\alpha,\nabla u_i\rangle\rangle-\rho\int_\Omega\langle\nabla
x_\alpha,\nabla u_i\rangle x_\alpha\Delta u_i\\
=&-\rho\int_\Omega\langle \nabla h_{\alpha i},\nabla \langle\nabla
x_\alpha,\nabla
u_i\rangle\rangle-\frac{\rho}{2}\int_\Omega\langle\nabla
x_\alpha^2,\nabla u_i\rangle \Delta u_i\\
\leq&(\delta\Lambda_i+n-2)\|\nabla h_{\alpha
i}\|^2+\frac{\rho^2}{4(\delta\Lambda_i+n-2)}\|\nabla \langle\nabla
x_\alpha,\nabla u_i\rangle\|^2\\
&-\frac{\rho}{2}\int_\Omega\langle\nabla x_\alpha^2,\nabla
u_i\rangle \Delta u_i.
\end{alignat}
 Applying
\eqref{2sec18} to \eqref{2sec17} yields
\begin{alignat}{1}\label{2sec19}
\sum_{i=1}^k(\Lambda_{k+1}-\Lambda_i)^2r_{\alpha
i}\leq&\sum_{i=1}^k(\Lambda_{k+1}-\Lambda_i)^2\Big(\delta p_{\alpha
i}+(\delta\Lambda_i+n-2)\|W_{\alpha
i}\|^2\\ \zdhy
&+(\delta-\rho)\|\langle\nabla x_\alpha,\nabla
u_i\rangle\|^2+\rho\|\langle\nabla x_\alpha,\nabla
u_i\rangle\|^2\Big)\\ \zdhy
&+\frac{1}{\delta}\sum_{i=1}^k(\Lambda_{k+1}-\Lambda_i)\|Z_{\alpha
i}\|^2\\ \zdhy
\leq&\sum_{i=1}^k(\Lambda_{k+1}-\Lambda_i)^2\Big(\delta p_{\alpha
i}+(\delta\Lambda_i+n-2)(\|W_{\alpha i}\|^2+\|\nabla h_{\alpha
i}\|^2)\\ \zdhy
&+(\delta-\rho)\|\langle\nabla x_\alpha,\nabla
u_i\rangle\|^2+\frac{\rho^2}{4(\delta\Lambda_i+n-2)}\|\nabla
\langle\nabla x_\alpha,\nabla
u_i\rangle\|^2\\ \zdhy
&-\frac{\rho}{2}\int_\Omega\langle\nabla x_\alpha^2,\nabla
u_i\rangle \Delta
u_i\Big)+\frac{1}{\delta}\sum_{i=1}^k(\Lambda_{k+1}-\Lambda_i)\|Z_{\alpha
i}\|^2\\ \zdhy
=&\sum_{i=1}^k(\Lambda_{k+1}-\Lambda_i)^2\Big(\delta p_{\alpha
i}+(\delta\Lambda_i+n-2)\|x_\alpha\nabla
u_i\|^2\\ \zdhy&
+(\delta-\rho)\|\langle\nabla x_\alpha,\nabla
u_i\rangle\|^2+\frac{\rho^2}{4(\delta\Lambda_i+n-2)}\|\nabla
\langle\nabla x_\alpha,\nabla
u_i\rangle\|^2\\ \zdhy
&-\frac{\rho}{2}\int_\Omega\langle\nabla x_\alpha^2,\nabla
u_i\rangle \Delta
u_i\Big)+\frac{1}{\delta}\sum_{i=1}^k(\Lambda_{k+1}-\Lambda_i)\|Z_{\alpha
i}\|^2.
\end{alignat}

Since $$\Delta h_{\alpha i}={\rm div}(\nabla h_{\alpha i})={\rm
div}(x_\alpha \nabla u_i)=\langle\nabla x_\alpha, \nabla
u_i\rangle+x_\alpha\Delta u_i,$$ we get from Proposition \ref{prop}
that $$\align \sum_{\alpha=1}^{n+1}p_{\alpha
i}=&\sum_{\alpha=1}^{n+1}\left(\int_{\Omega}h_{\alpha
i}(-\Delta)^ph_{\alpha i}-\Lambda_i\|x_\alpha\nabla
u_i\|^2-\|\langle\nabla x_\alpha,\nabla u_i\rangle\|^2\right)\\
=&\sum_{\alpha=1}^{n+1}\int_{\Omega}h_{\alpha
i}(-\Delta)^ph_{\alpha i}-(\Lambda_i+1)\\
=&\sum_{\alpha=1}^{n+1}\int_{\Omega}\left(\langle\nabla
x_\alpha, \nabla u_i\rangle+x_\alpha\Delta
u_i\right)(-\Delta)^{p-2}\left(\langle\nabla x_\alpha, \nabla
u_i\rangle+x_\alpha\Delta u_i\right)-(\Lambda_i+1)\\
\leq&f(\Lambda_i, n)-(\Lambda_i+1).
\endalign$$
A direct calculation yields (see (2.44), (2.45), (2.46) and (2.47)
in \cite{wangxia07})
$$\align
&\sum_{\alpha=1}^{n+1}r_{\alpha i}=n,
\\
&\sum_{\alpha=1}^{n+1}\|x_\alpha\nabla
u_i\|^2=\sum_{\alpha=1}^{n+1}\|\langle\nabla x_\alpha,\nabla
u_i\rangle\|^2=1,
\\
&\sum_{\alpha=1}^{n+1}\|\nabla \langle\nabla x_\alpha,\nabla
u_i\rangle\|^2=\Lambda_i-(n-2),
\\
&\sum_{\alpha=1}^{n+1}\|Z_{\alpha i}\|^2=\Lambda_i+\frac{(n-2)^2}{
4}.\endalign
$$
Therefore, summing up \eqref{2sec19} over $\alpha$ from 1 to
$n+1$, one gets
$$\align
&n\sum_{i=1}^k(\Lambda_{k+1}-\Lambda_i)^2\\
\leq&\sum_{i=1}^k(\Lambda_{k+1}-\Lambda_i)^2
\Big(\delta\left(f(\Lambda_i, n)-(\Lambda_i+1)\right)+(\delta\Lambda_i+n-2)+(\delta-\rho)\\
&+\frac{\rho^2}{4(\delta\Lambda_i+n-2)}\left(\Lambda_i-(n-2)\right)\Big)
+\frac1\delta\sum_{i=1}^k(\Lambda_{k+1}-\Lambda_i)\left(\Lambda_i+\frac{(n-2)^2}{
4}\right).
\endalign$$
That is,
\begin{alignat}{1}\label{2sec20}
&2\sum_{i=1}^k(\Lambda_{k+1}-\Lambda_i)^2\\
\leq&\sum_{i=1}^k(\Lambda_{k+1}-\Lambda_i)^2
\left(\delta f(\Lambda_i,
n)-\rho+\frac{\rho^2}{4(\delta\Lambda_i+n-2)}\left(\Lambda_i-(n-2)\right)\right)\\
&+\frac1\delta\sum_{i=1}^k(\Lambda_{k+1}-\Lambda_i)\left(\Lambda_i+\frac{(n-2)^2}{
4}\right).
\end{alignat}
Taking $$\rho=\frac{2(\delta\Lambda_i+n-2)}{\Lambda_i-(n-2)}$$ in
\eqref{2sec20} yields
$$\align
2\sum_{i=1}^k(\Lambda_{k+1}-\Lambda_i)^2\leq&\sum_{i=1}^k(\Lambda_{k+1}-\Lambda_i)^2
\left(\delta f(\Lambda_i,
n)-\frac{\delta\Lambda_i+n-2}{\Lambda_i-(n-2)}\right)\\
&+\frac1\delta\sum_{i=1}^k(\Lambda_{k+1}-\Lambda_i)\left(\Lambda_i+\frac{(n-2)^2}{
4}\right).
\endalign$$ Hence, we obtain
\begin{alignat}{1}\label{2sec21}
&\sum_{i=1}^k(\Lambda_{k+1}-\Lambda_i)^2
\left(2+\frac{n-2}{\Lambda_i-(n-2)}\right)\\
\leq
&\delta\sum_{i=1}^k(\Lambda_{k+1}-\Lambda_i)^2 \left(f(\Lambda_i,
n)-\frac{\Lambda_i}{\Lambda_i-(n-2)}\right)\\
&+\frac{1}{\delta}\sum_{i=1}^k(\Lambda_{k+1}-\Lambda_i)
\left(\Lambda_i+\frac{(n-2)^2}{4}\right).
\end{alignat}
Minimizing the right hand side of \eqref{2sec21} as a function of
$\delta$ by choosing
$$\delta=\left(\frac{\sum\limits_{i=1}^k(\Lambda_{k+1}-\Lambda_i)
\left(\Lambda_i+\frac{(n-2)^2}{4}\right)}{\sum\limits_{i=1}^k(\Lambda_{k+1}-\Lambda_i)^2
\left(f(\Lambda_i, n)-\frac{\Lambda_i}{\Lambda_i-(n-2)}\right)}\right)^{\frac12}$$
concludes the proof of Theorem \ref{thm}.

{\bf Proof of Corollary \ref{coro}.}

It is easy to see from \eqref{Intr9} that
\begin{alignat}{1}\label{2sec22}
\sum_{i=1}^k(\Lambda_{k+1}-\Lambda_i)^2
\leq&\left\{\sum_{i=1}^k(\Lambda_{k+1}-\Lambda_i)^2
\left(f(\Lambda_i,
n)-\frac{\Lambda_i}{\Lambda_i-(n-2)}\right)\right\}^{\frac{1}{2}}\\
&\hspace{2.6cm}\times\left\{\sum_{i=1}^k(\Lambda_{k+1}-\Lambda_i)
\left(\Lambda_i+\frac{(n-2)^2}{4}\right)\right\}^{\frac{1}{2}}.
\end{alignat}
One can check by induction that
$$\align
&\left\{\sum_{i=1}^k(\Lambda_{k+1}-\Lambda_i)^2 \left(f(\Lambda_i,
n)-\frac{\Lambda_i}{\Lambda_i-(n-2)}\right)\right\}
\left\{\sum_{i=1}^k(\Lambda_{k+1}-\Lambda_i)
\left(\Lambda_i+\frac{(n-2)^2}{4}\right)\right\}\\
\leq&\left(\sum_{i=1}^k(\Lambda_{k+1}-\Lambda_i)^2\right)
\left\{\sum_{i=1}^k(\Lambda_{k+1}-\Lambda_i)\left(f(\Lambda_i,
n)-\frac{\Lambda_i}{\Lambda_i-(n-2)}\right)
\left(\Lambda_i+\frac{(n-2)^2}{4}\right)\right\},
\endalign$$
which together with \eqref{2sec22} yields inequality \eqref{Intr10}.

Solving the quadratic polynomial of $\Lambda_{k+1}$ in
\eqref{Intr10}, we obtain inequality \eqref{Intr11} and
\eqref{Intr12}. It completes the proof of Corollary \ref{coro}.

\end{document}